\documentclass[fleqn]{article}
\usepackage{graphicx}
\usepackage{amssymb,amsmath}
\usepackage{theorem}
\usepackage{hyperref}

\newtheorem{theorem}{Theorem}[section]
\newtheorem{prop}[theorem]{Proposition}
\newtheorem{lemma}[theorem]{Lemma}
\newtheorem{claim}[theorem]{Claim}
\newtheorem{corollary}[theorem]{Corollary}

\theorembodyfont{\rmfamily}

\newenvironment{proof}{\noindent{\bf Proof. }}{\hfill$\square$\medskip}
\begin{document}

\addtolength{\baselineskip}{3pt} \setlength{\oddsidemargin}{0.2in}

\def\R{{\mathbb R}}
\def\one{{\mathbf 1}}
\def\Q{{\mathbb Q}}
\def\Z{{\mathbb Z}}
\def\C{{\mathbb C}}
\def\N{{\mathbb N}}
\def\hom{{\rm hom}}
\def\Hom{{\rm Hom}}
\def\Inj{{\rm Inj}}
\def\Ind{{\rm Ind}}
\def\inj{{\rm inj}}
\def\ind{{\rm ind}}
\def\sur{{\rm sur}}
\def\PAG{{\rm PAG}}
\def\CUT{{\text{\rm CUT}}}
\def\eps{\varepsilon}
\long\def\killtext#1{}
\def\dd{{d_{\text{\rm set}}}}
\def\dl{{d_{\text{\rm left}}}}
\def\dr{{d_{\text{\rm right}}}}
\def\Ge{{\mathbf G}}

\def\ontop#1#2{\genfrac{}{}{0pt}{}{#1}{#2}}

\def\HH{{\cal H}}
\def\SS{{\cal S}}
\def\TT{{\cal T}}
\def\II{{\cal I}}
\def\FF{{\cal F}}
\def\GG{{\cal G}}
\def\AA{{\cal A}}
\def\MM{{\cal M}}
\def\EE{{\cal E}}
\def\LL{{\cal L}}
\def\PP{{\cal P}}
\def\CC{{\cal C}}
\def\SS{{\cal S}}
\def\QG{{\cal QG}}
\def\QGM{{\cal QGM}}
\def\CD{{\cal CD}}
\def\P{{\sf P}}
\def\E{{\sf E}}
\def\Var{{\sf Var}}
\def\T{^{\sf T}}

\def\tr{{\rm tr}}
\def\cost{\hbox{\rm cost}}
\def\val{\hbox{\rm val}}
\def\rk{{\rm rk}}
\def\diam{{\rm diam}}
\def\Ker{{\rm Ker}}

\title{\huge\bf Limits of dense graph sequences\\[12mm]}

\author{
{\sc L\'aszl\'o Lov\'asz} and
{\sc Bal\'azs Szegedy}\\[8mm]
Microsoft Research \\
One Microsoft Way\\
Redmond, WA 98052}

\date{Technical Report TR-2004-79\\
August 2004}

\maketitle

\tableofcontents

\begin{abstract}
We show that if a sequence of dense graphs $G_n$ has the property
that for every fixed graph $F$, the density of copies of $F$ in $G_n$
tends to a limit, then there is a natural ``limit object'', namely a
symmetric measurable function $W:~[0,1]^2\to[0,1]$. This limit object
determines all the limits of subgraph densities. Conversely, every
such function arises as a limit object. We also characterize graph
parameters that are obtained as limits of subgraph densities by the
``reflection positivity'' property.

Along the lines we introduce a rather general model of random graphs,
which seems to be interesting on its own right.
\end{abstract}

\section{Introduction}

Let $G_n$ be a sequence of simple graphs whose number of nodes tends
to infinity. For every fixed simple graph $F$, let $\hom(F,G)$ denote
the number homomorphisms of $F$ into $G$ (edge-preserving maps
$V(F)\to V(G)$). We normalize this number to get the {\it
homomorphism density}
\[
t(F,G)=\frac{\hom(F,G)}{|V(G)|^{|V(F)|}}.
\]
This quantity is the probability that a random mapping $V(F)\to V(G)$
is a homomorphism.

Suppose that the graphs $G_n$ become more and more similar in the
sense that $t(F,G_n)$ tends to a limit $t(F)$ for every $F$. Let
$\TT$ denote the set of graph parameters $t(F)$ arising this way. The
goal of this paper is to give characterizations of graph parameters
in $\TT$. (This question is only interesting if the graphs $G_n$ are
dense (i.e., they have $\Omega(|V(G_n)|^2)$ edges); else, the limit
is 0 for every simple graph $F$ with at least one edge.)

One way to characterize members of $\TT$ is to define an appropriate
limit object from which the values $t(F)$ can be read off.

For example, let $G_n$ be a random graph with density $1/2$ on $n$
nodes. It can be shown that this converges with probability 1. A
natural guess for the limit object would be the random countable
graph. This is a very nice object, uniquely determined up to
automorphism. However, this graph is too ``robust'': the limit of
random graphs with edge-density $1/3$ would be the same, while the
homomorphism densities have different limits than in the case of
edge-density $1/2$.

The main result of this paper is to show that indeed there is a
natural ``limit object'' in the form of a symmetric measurable
function $W:~[0,1]^2\to[0,1]$ (we call $W$ {\it symmetric} if
$W(x,y)=W(y,x)$). Conversely, every such function arises as the limit
of an appropriate graph sequence. This limit object determines all
the limits of subgraph densities: if $F$ is a simple graph with
$V(F)=[k]=\{1,\dots,k\}$, then
\begin{equation}\label{INTREP}
t(f,W)=\int\limits_{[0,1]^n} \prod_{ij\in E(F)} W(x_1,x_j)\,dx_1\dots
dx_n.
\end{equation}
The limit object for random graphs of density $p$ is the constant
function $p$.

Another characterization of graph parameters $t(F)$ that are limits
of homomorphism densities can be given by describing a complete
system of inequalities between the values $t(F)$ for different finite
graphs $F$. One can give such a complete system in terms of the
positive semidefiniteness of a certain sequence of matrices which we
call {\it connection matrices} (see section \ref{REFPOS} for
details). This property is related to {\it reflection positivity} in
statistical mechanics. Our results in this direction can be thought
of as analogues of the characterization of homomorphism density
functions given in \cite{FLS} in the limiting case.

We can also look at this result as an analogue of the well known
characterization of moment sequences in terms of the positive
semidefiniteness of the moment matrix. A ``2-variable'' version of a
sequence is a graph parameter, and representation in form of moments
of a function (or random variable) can be replaced by the integral
representation (\ref{INTREP}). The positive semidefiniteness of
connection matrices is analogous to the positive semidefiniteness of
moment matrices.

\medskip

Every symmetric measurable function $W:~[0,1]^2\to[0,1]$ gives rise
to a rather general model of random graphs, which we call $W$-random.
Their main role in this paper is that they provide a graph sequence
that converges to $W$; but they seem to be interesting on their own
right.

We show that every random graph model satisfying some rather natural
criteria can be obtained as a $W$-random graph for an appropriate
$W$.

\medskip

The set $\TT$ was introduced by Erd\H{o}s, Lov\'asz and Spencer
\cite{ELS}, where the dimension of its projection to any finite
number of coordinates (graphs $F$) was determined.

Limit objects of graph sequences were constructed by Benjamini and
Schramm \cite{BS} for sequences of graphs with bounded degree;
this was extended by Lyons \cite{Ly} to sequences of graphs with
bounded average degree. (The normalization in that case is
different.)

\section{Definitions and main results}\label{RESULTS}

\subsection{Weighted graphs and homomorphisms}

A {\it weighted graph} $G$ is a graph with a weight $\alpha_G(i)$
associated with each node and a weight $\beta_G(i,j)$ associated with
each edge $ij$. (Here we allow that $G$ has loops, but no multiple
edges.) In this paper we restrict our attention to positive real
weights between 0 and 1. An edge with weight 0 will play the same
role as no edge between those nodes, so we could assume that we only
consider weighted complete graphs with loops at all nodes (but this
is not always convenient). The adjacency matrix of a weighted graph
is obtained by replacing the 1's in the adjacency matrix by the
weights of the edges. An {\it unweighted graph} is a weighted graph
where all the node- and edgeweights are 1. We set
\[
\alpha_G=\sum_{i\in V(G)} \alpha_G(i).
\]

Recall that for two unweighted graphs $F$ and $G$, $\hom(F,G)$
denotes the number of homomorphisms (adjacency preserving maps) from
$F$ to $G$. We extend this notion to the case when $G$ is a weighted
graph. To every $\phi:~V(F)\to V(G)$, we assign the weights
\[
\alpha_\phi=\prod_{u\in V(F)} \alpha_G(\phi(u)),
\]
and
\[
\hom_{\phi}(F,G)=\prod_{uv\in E(F)} \beta_G(\phi(u),\phi(v)).
\]
We then define the {\it homomorphism function}
\begin{equation}\label{HOMSUM}
\hom(F,G)=\sum_{\phi:~V(F)\to V(G)} \alpha_\phi \hom_\phi(F,G).
\end{equation}
and the {\it homomorphism density}
\[
t(F,G)=\frac{\hom(F,G)}{\alpha_G^{|V(F)|}}.
\]
We can also think of $t(F,G)$ as a homomorphism function after the
nodeweights of $G$ are scaled so that their sum is 1.

It will be convenient to extend the notation $\hom_\phi$ as follows.
Let $\phi:~V'\to V(G)$ be a map defined on a subset $V'\subseteq
V(F)$. Then define
\[
\alpha_\phi=\prod_{u\in V'} \alpha_G(\phi(u)),
\]
and
\[
\hom_\phi(F,G)= \sum_{\ontop{\psi: V(F)\to
V(G)}{\psi\text{~extends~}\phi}}
\frac{\alpha_\psi}{\alpha_\phi}\hom_\psi(F,G).
\]
If $V'=\emptyset$, then $\alpha_\phi=1$ and
$\hom_\phi(F,G)=\hom(F,G)$.

\subsection{Convergence of graph sequences}

Let $(G_n)$ be a sequence of weighted graphs. We say that this
sequence is {\it convergent}, if the sequence $(t(F,G_n))$ has a
limit as $n\to\infty$ for every simple unweighted graph $F$. (Note
that it would be enough to assume this for connected graphs $F$.) We
say that the sequence {\it converges to a finite weighted graph $G$}
if
\[
t(F,G_n)\longrightarrow \hom(F,G)
\]
for every simple graph $F$. A convergent graph sequence may not
converge to any finite weighted graph; it will be our goal to
construct appropriate limit objects for convergent graph sequences
which do not have a finite graph as a limit.

\medskip

A {\it graph parameter} is a function defined on simple graphs that
is invariant under isomorphism. Every weighted graph $G$ defines
graph parameters $\hom(.,G)$, $\inj(.,G)$, $t(.,G)$ and $t_0(.,G)$.

Often we can restrict our attention to graph parameters $f$
satisfying $f(K_1)=1$, which we call {\it normalized}. Of the four
parameters above, $t(.,G)$ and $t_0(.,G)$ are normalized. We say that
a graph parameter $f$ is {\it multiplicative}, if
$f(G_1G_2)=f(G_1)f(G_2)$, where $G_1G_2$ denotes the disjoint union
of two graphs $G_1$ and $G_2$. The parameters $\hom(.,G)$ and
$t(.,G)$ are multiplicative.

The same graph parameter $\hom(.,G)$, defined by a weighted graph
$G$, arises from infinitely many graphs. Replace a node $i$ of $G$ by
two nodes $i'$ and $i''$, whose weights are chosen so that
$\alpha(i')+\alpha(i'')=\alpha(i)$; define the edge weights
$\beta(i'j)=\beta(i''j)=\beta(ij)$ for every node $j$; and keep all
the other nodeweights and edgeweights. The resulting weighted graph
$G'$ will define the same graph parameter $\hom(.,G')=\hom(.,G)$.
Repeating this operation we can create arbitrarily large weighted
graphs defining the same graph parameter.

If we want to stay among unweighted graphs, then the above operation
cannot be carried out, and the function $\hom(.,G)$ in fact
determines $G$ \cite{LL}. But for the $t(.,G)$ parameter, the
situation is different: if we replace each node of an unweighted
graph $G$ by $N$ copies (where copies of two nodes are connected if
and only if the originals were), then $t(.,G)$ does not change.

In particular, if we consider a convergent graph sequence, we need
not assume that the number of nodes tends to infinity: we could
always achieve this without changing the limit.

\subsection{Reflection positivity}\label{REFPOS}

A {\it $k$-labeled graph} ($k\ge 0$) is a finite graph in which $k$
nodes are labeled by $1,2,\dots k$ (the graph can have any number of
unlabeled nodes). For two $k$-labeled graphs $F_1$ and $F_2$, we
define the graph $F_1F_2$ by taking their disjoint union, and then
identifying nodes with the same label, and then cancelling the
resulting multiplicities of edges. Hence for two $0$-labeled graphs,
$F_1F_2$ is just their disjoint union.

Let $f$ be any graph parameter defined on simple graphs. For every
integer $k\ge 0$, we define the {\it connection matrix} $M(k,f)$ as
follows. This is infinite matrix, whose rows and columns are indexed
by (isomorphism types of) $k$-labeled graphs. The entry in the
intersection of the row corresponding to $F_1$ and the column
corresponding to $F_2$ is $f(F_1F_2)$. We say that the parameter $f$
is {\it reflection positive}, if $M(k,f)$ is positive semidefinite
for every $k\ge 0$.

In \cite{FLS}, a related matrix was defined. In that paper, the test
graphs $F$ may have multiple edges and the target graphs $G$ have
arbitrary edgeweights. Let us call a graph parameter defined on
graphs which may have multiple edges a {\it multigraph parameter}.
The only difference in the definition of the connection matrix is
that edge multiplicities are not cancelled when $F_1F_2$ is defined.
It was shown that $\hom(.,G)$ as a multigraph parameter is reflection
positive for every weighted graph $G$, and the matrix $M(f,k)$ has
rank at most $|V(G)|^k$. It was also shown that these two properties
characterize which multigraph parameters arise in this form. (In this
paper we restrict our attention to simple test graphs $F$ and
edge-weights between 0 and 1. See also section \ref{EXTENSIONS}.)

For graph parameters defined on simple graphs, there is a simpler
matrix whose positive semidefiniteness could be used to define
reflection positivity. Let $M_0(k,f)$ denote the submatrix of
$M(k,f)$ formed by those rows and columns that are indexed by
$k$-labeled graphs on $k$ nodes (so that every node is labeled). The
equivalence of these definitions (under some further conditions) will
follow from our main theorem.

We could combine all these matrices into single matrix $M_0(f)$: the
rows and columns of $M_0(f)$ are indexed by all finite graphs whose
nodes form a finite subset of $\N$. To get the entry in the
intersection of row $F_1$ and column $F_2$, we take the union
$F_1\cup F_2$, and evaluate $f$ on this union. Clearly every
$M_0(k,f)$ is a minor of $M_0(f)$, and every finite minor of $M_0(f)$
is a minor of $M_0(k,f)$ for every large enough $k$. So $M_0(f)$ is
positive semidefinite if and only if every $M_0(k,f)$ is.

\subsection{Homomorphisms, subgraphs, induced subgraphs}

Sometimes it is more convenient to work with injective maps. For two
unweighted graphs $F$ and $G$, let $\inj(F,G)$ denote the number of
injective homomorphisms from $F$ to $G$ (informally, the number of
copies of $F$ in $G$). We also introduce the {\it injective
homomorphism density}
\[
t_0(F,G)=\frac{\inj(F,G)}{(|V(G)|)_{|V(F)|}}
\]
(where $(n)_k=n(n-1)\cdot(n-k+1)$).

From a graph-theoretic point of view, it is also important to count
induced subgraphs. More precisely, if $F$ and $G$ are two unweighted
graphs, then let $\ind(F,G)$ denote the number of embeddings of $F$
into $G$ as an induced subgraph. We define the {\it induced
homomorphism density} by
\[
t_1(F,G)=\frac{\ind(F,G)}{(|V(G)|)_{|V(F)|}}.
\]

If $G$ is weighted, then we define $\inj(F,G)$ by the same type of
sum as for homomorphisms:
\[
\inj(F,G)=\sum_{\phi} \alpha_\phi \hom_\phi(F,G),
\]
except that the summation is restricted to injective maps. Let
\[
t_0(F,G)=\frac{\inj(F,G)}{(\alpha_G)_{|V(F)|}}
\]
(where for a sequence $\alpha=(\alpha_1,\dots,\alpha_n)$,
$(\alpha)_k$ denotes the $k$-th elementary symmetric polynomial of
the $\alpha_i$). Note that the normalization was chosen so that
$t(F,G)=t_0(F,G)=1$ if $F$ has no edges.

We also extend the $\ind$ function to the case when $G$ is weighted:
we define
\[
\ind_{\phi}(F,G)=\prod_{uv\in E(F)} \beta_G(\phi(u),\phi(v))
\prod_{uv\in E(\overline{F})} (1-\beta_G(\phi(u),\phi(v)))
\]
(here $\overline{F}$ denotes the complement of the graph $F$),
\[
\ind(F,G)=\sum_{\phi} \alpha_\phi \ind_\phi(F,G),
\]
and
\[
t_1(F,G)=\frac{\ind(F,G)}{(\alpha_G)_{|V(F)|}}.
\]

In the definition of convergence, we could replace $t(F,G)$ by the
number $t_0(F,G)$ of embeddings (injective homomorphisms); this is
more natural from the graph theoretic point if view. This would not
change the notion of convergence or the value of the limit, as the
following simple lemma shows:

\begin{lemma}\label{HOM-INJ}
For every weighted graph $G$ and unweighted simple graph $F$, we have
\[
|t(F,G) -t_0(F,G)| <\frac{1}{|V(G)|}\binom{|V(F)|}{2}.
\]
\end{lemma}

We could also replace the $\hom$ function by $\ind$ function. Indeed,
\[
\inj(F,G)=\sum_{F'\supset F}\ind(F',G)
\]
(where $F'$ ranges over all supergraphs  of $F$ on the same set of
nodes), and by inclusion-exclusion,
\[
\ind(F,G)=\sum_{F'\supset F}(-1)^{|E(F')\setminus E(F)|} \ind(F',G).
\]
Hence it follows that
\[
t_1(F,G)=\sum_{F'\supset F}(-1)^{|E(F')\setminus E(F)|} t_0(F',G).
\]
It will be convenient to introduce the following operator: if $f$ is
any graph parameter, then $f^{\dagger}$ is the graph parameter
defined by
\[
f^{\dagger}(F)=\sum_{F'\supset F}(-1)^{|E(F')\setminus E(F)|} f(F').
\]
Thus $t_1=t_0^\dagger$.

(There is a similar precise relation between the numbers of
homomorphisms and injective homomorphisms as well, but we will not
have to appeal to it.)

\subsection{The limit object}

We'll show that every convergent graph sequence has a limit object,
which can be viewed as an infinite weighted graph on the points of
the unit interval. To be more precise, for every symmetric measurable
function $W:~[0,1]^2\to[0,1]$, we can define a graph parameter
$t(.,W)$ by
\[
t(F,W)=\int_{[0,1]^k} \prod_{ij\in E(F)}
W(x_i,x_j)\,dx_1\,\dots\,dx_k
\]
(where $F$ is a simple graph with $V(F)=[k]$).

Let $\phi:~V'\to [0,1]$ is a map defined on a subset $V'\subseteq
[k]$. Similarly as in the case of homomorphisms into finite graphs,
we define $t_\phi$ as follows. Let, say, $V'=[k']$ ($1\le k'\le k$),
and $x_i=\phi(i)$ ($i=1,\dots,k'$). Define
\[
t_\phi(F,W)=\int_{[0,1]^{k-k'}} \prod_{ij\in E(F)}
W(x_i,x_j)\,dx_{k'+1}\,\dots\,dx_k.
\]

It is easy to see that for every weighted finite graph $H$, the
simple graph parameter $t(.,H)$ is a special case. Indeed, define a
function $W_H:~[0,1]^2\to[0,1]$ as follows. Let $\alpha_i$ be the
nodeweights of $H$ and $\beta_{ij}$, the edgeweights of $H$. We may
assume that $\sum_i\alpha_i=1$. For $(x,y)\in[0,1]^2$, let $a$ and
$b$ determined by
\begin{align*}
\alpha_1+\dots+\alpha_{a-1} &\le x <\alpha_1+\dots+\alpha_a,\\
\alpha_1+\dots+\alpha_{b-1} &\le y <\alpha_1+\dots+\alpha_b,
\end{align*}
and let
\[
W_H(x,y)=\beta_{ab}.
\]

The main result in this paper is the following. Recall that $\TT$
denotes the set of graph parameters $f$ that are limits of graph
parameters $t(.,G)$; i.e., there is a convergent sequence of simple
graphs $G_n$ such that
\[
f(.)=\lim_{n\to\infty} t(.,G_n)
\]
for every simple graph $F$.

\begin{theorem}\label{MAIN}
For a simple graph parameter $f$ the following are equivalent:

\smallskip

{\rm (a)} $f\in\TT$;

\smallskip

{\rm (b)} There is a symmetric measurable function
$W:~[0,1]^2\to[0,1]$ for which $f=t(.,W)$.

\smallskip

{\rm (c)} The parameter $f$ is normalized, multiplicative and
reflection positive.

\smallskip

{\rm (d)} The parameter $f$ is normalized, multiplicative and
$M_0(f)$ is positive semidefinite.

\smallskip

{\rm (e)} The parameter $f$ is normalized, multiplicative and
$f^{\dagger}\ge 0$.
\end{theorem}

\medskip

\noindent{\bf Remarks 1.} The theorem gives four characterizations of
the set $\TT$: one analytic, two algebraic and one combinatorial.
Characterizations (c), (d) and (e) are closely related (even though a
direct proof of the equivalence of (c) and (d) is not easy). Any of
these three on the one hand, and (b) on the other, form a ``dual''
pair in the spirit of NP-coNP: (b) tells us why a graph parameter is
in $\TT$, while (c) (or (d) or (e)) tells us why it is not.

\medskip

{\bf 2.} Corollary \ref{W-CONVERGE} below shows that finite weighted
graphs are limits of simple unweighted graphs. This implies that in
the definition of $\TT$, we could take convergent sequences of
weighted graphs instead of unweighted graphs.

\medskip

{\bf 3.} We could define a more general limit object as a probability
space $(\Omega,\AA,\pi)$ and a symmetric measurable function on
$W:~\Omega\times\Omega\to[0,1]$. This would not give rise to any new
invariants. However, some limit objects may have a simpler or more
natural representation on other $\sigma$-algebras (cf. Corollary
\ref{W-MULT} below).

\medskip

{\bf 4.} One might think that (c) and (d) are equivalent for the
more direct reason that $M_0(k,f)$ is positive semidefinite if and
only if $M(k,f)$ is. This implication, however, does not hold for
a fixed $k$ (even if we assume that $f$ is normalized and
multiplicative). For example, $M_0(1,f)$ is positive semidefinite
for every normalized multiplicative graph parameter, but $M(1,f)$
is not if $f(F)$ is the number of matchings in $F$.

\medskip

{\bf 5.} In the case when $f=\hom(F,G)$ for some finite graph $G$,
$f^{\dagger}\ge 0$ in condition (e) expresses that counting induced
subgraphs in $G$ we get non-negative values.

\medskip

As an immediate application of Theorem \ref{MAIN}, we prove the
following fact:

\begin{prop}\label{W-MULT}
If $t_1,t_2\in \TT$, then $t_1t_2\in\TT$.
\end{prop}

This follows from condition (c) in Theorem \ref{MAIN}, using that
positive semidefiniteness is preserved under Schur product. It may be
instructive to see how a representation of the product of type (b)
can be constructed. Let $t_i=t(.,W_i)$, and define $W$ as the
4-variable function $W_1(x_1y_1)W_2(x_2,y_2)$. We can consider $W$ as
a function in two variables $x,y$, where $x=(x_1,x_2)\in[0,1]^2$ and
$y=(y_1,y_2)\in[0,1]^2$. Then $W$ gives rise to graph parameter
$t(.,W)$, and it is straightforward to check that $t=t_1t_2$.

\subsection{$W$-random graphs}

Given any symmetric measurable function $W:~[0,1]^2\to[0,1]$ and an
integer $n>0$, we can generate a random graph $\Ge(n,W)$ on node set
$[n]$ as follows. We generate $n$ independent numbers $X_1,\dots,X_n$
from the uniform distribution on $[0,1]$, and then connect nodes $i$
and $j$ with probability $W(X_i,X_j)$.

As a special case, if $W$ is the identically $p$ function, we get
``ordinary'' random graphs $\Ge(n,p)$. This sequence is convergent
with probability 1, and in fact it converges to the graph $K_{1}(p)$,
the weighted graph with one node and one loop with weight $p$. The
limiting simple graph parameter is given by $t(F)=p^{|E(F)|}$.

More generally, let $W_H:~[0,1]^2\to[0,1]$ be defined by a (finite)
weighted graph $H$ with $V(H)=[q]$, whose node weights $\alpha_i$
satisfy $\alpha_1+\dots+\alpha_q=1$. Then $\Ge(n,W_H)$ can be
described as follows. We open $q$ bins $V_1,\dots, V_q$. Create $n$
nodes, and put each of them independently in bin $i$ with probability
$\alpha_i$. For every pair $u,v$ of nodes, connect them by an edge
with probability $\beta_{ij}$ if $u\in V_i$ and $v\in V_j$. We call
$\Ge(n,W_H)$ a {\it random graph with model} $H$.

We show that the homomorphisms densities into $\Ge(n,F)$ are close to
the homomorphism densities into $W$. Let us fix a simple graph $F$,
let $V(F)=[k]$ and $\Ge=\Ge(n,W)$.

The following lemma summarizes some simple properties of $W$-random
graphs.

\begin{lemma}\label{EXPINJ}
For every simple graph $F$,

\smallskip

{\rm (a)} $\E\bigl(t_0(F,\Ge(n,W))\bigr)= t(F,W)$;

\smallskip

{\rm (b)} $\displaystyle \bigl|\E\bigl(t(F,\Ge(n,W))\bigr)-
t(F,W)\bigr|<\frac{1}{n}\binom{|V(F)|}{2}$;

\smallskip

{\rm (c)} $\displaystyle \Var(t(F,\Ge(n,W)))\le\frac{3}{n}|V(F)|^2$.
\end{lemma}

This lemma implies, by Chebyshev's inequality, that
\[
\Pr\bigl(|t(F,\Ge(n,W))-t(F,W)|>\eps\bigr) \leq
3|V(F)|^2\frac{1}{n\eps^2}.
\]
Much stronger concentration results can be proved for $t(F,\Ge)$,
using deeper techniques (Azuma's inequality):

\begin{theorem}\label{CONC}
Let $F$ be a graph with $k$ nodes. Then for every $0<\eps<1$,
\begin{equation}\label{INJFH}
\Pr\Bigl(|t_0(F,\Ge(n,W))- t(F,W)| > \eps\Bigr) \le
2\exp\left(-\frac{\eps^2}{2k^2}n\right).
\end{equation}
and
\begin{equation}\label{TFH}
\Pr\Bigl(|t(F,\Ge(n,W))-t(F,W)| > \eps \Bigr) \le
2\exp\left(-\frac{\eps^2}{18k^2}n\right).
\end{equation}
\end{theorem}

From this Theorem it is easy to show:

\begin{corollary}\label{W-CONVERGE}
The graph sequence $\Ge(n,W)$ is convergent with probability 1, and
its limit is the function $W$.
\end{corollary}

Indeed, the sum of the right hand sides is convergent for every fixed
$\eps>0$, so it follows by the Borell--Cantelli Lemma that
$t(F,\Ge(n,W))\to t(F,W)$ with probability 1. There is only a
countable number of graphs $F$, so this holds with probability 1 for
every $F$.

This way of generating random graphs is quite general in the
following sense. Suppose that for every $n\ge 1$, we are given a
distribution on simple graphs on $n$ given nodes, say $[n]$; in other
words, we have a random variable $\Ge_n$ whose values are simple
graphs on $[n]$. We call this random variable a {\it random graph
model}. Clearly, every symmetric function $W:~[0,1]^2\to[0,1]$ gives
rise to a random graph model $\Ge(n,W)$. The following theorem shows
that every model satisfying rather general conditions is of this
form:

\begin{theorem}\label{GEN-RAND}
A random graph model is of the form $\Ge(n,W)$ for some symmetric
function $W:~[0,1]^2\to[0,1]$ if and only if it has the following
three properties:

\smallskip

{\rm (i)} the distribution of $\Ge_n$ is invariant under relabeling
nodes;

\smallskip

{\rm (ii)} if we delete node $n$ from $\Ge_n$, the distribution of
the resulting graph is the same as the distribution of $\Ge_{n-1}$;

\smallskip

{\rm (iii)} for every $1<k<n$, the subgraphs of $\Ge$ induced by
$[k]$ and $\{k+1,\dots,n\}$ are independent as random variables.
\end{theorem}

\section{Examples}

\subsection{Quasirandom graphs}\label{RANDOM}

Graph sequences converging to $K_1(p)$ are well studied under the
name of {\it quasirandom graphs} (see \cite{CGW}). More generally,
graph sequences converging to a finite weighted graph $H$ are called
{\it quasirandom graphs with model} $H$. The name is again justified
since random graphs $\Ge(n,W_H)$ with model $H$ converge to $H$ with
probability 1. These generalized quasirandom graphs are characterized
in \cite{LS}.

\subsection{Half-graphs}\label{HALF}

Let $H_{n,n}$ denote the bipartite graph on $2n$ nodes
$\{1,\dots,n,1',\dots,n'\}$, where $i$ is connected to $j'$ if and
only if $i\le j$. It is easy to see that this sequence is convergent.
Indeed, let $F$ be a simple graph with $k$ nodes; we show that the
limit of $t(F,H_{n,n})$ exists. We may assume that $F$ is connected.
If $F$ is non-bipartite, then $t(F,H_{n,n})=0$ for all $n$, so
suppose that $F$ is bipartite; let $V(F)=V_1\cup V_2$ be its (unique)
bipartition. Then every homomorphisms of $F$ into $H$ preserves the
2-coloring, and so the homomorphisms split into two classes: those
that map $V_1$ into $\{1,\dots,n\}$ and those that map it into
$\{1',\dots,n'\}$. By the symmetry of the half-graphs, these two
classes have the same cardinality.

Now $F$ defines a partial order $P$ on $V(F)$, where $u\le v$ if and
only if $u=v$ or $u\in V_1$, $v\in V_2$, and $uv\in E$.
$(1/2)\hom(F,H_{n,n})$ is just the number of order-preserving maps of
$P$ to the chain $\{1,\dots,n\}$, and so
\[
\frac{1}{2^{k-1}} t(F,H_{n,n}) = \frac{1}{2^{k-1}} \cdot
\frac{\hom(F,H_{n,n})}{(2n)^k}=\frac{(1/2)\hom(F,H_{n,n})}{n^k}
\]
is the probability that a random map of $V(F)$ into $\{1,\dots,n\}$
is order-preserving. As $n\to\infty$, the fraction of non-injective
maps tends to $0$, and hence it is easy to see that
$2^{1-k}t(F,H_{n,n})$ tends to a number $2^{1-k}t(F)$, which is the
probability that a random ordering of $V(F)$ is compatible with $P$.
In other words, $k!2^{1-k}t(F)$ is the number of linear extensions of
$P$.

However, the half-graphs do not converge to any finite weighted
graph. To see this, let $S_k$ denote the star on $k$ nodes, and
consider the (infinite) matrix $M$ defined $M_{k,l}=t(S_{k+l-1})$. If
$t(F)=t(F,G_0)$ for some finite weighted graph $G_0$, then it follows
from the characterization of homomorphism functions in \cite{FLS}
that this matrix has rank at most $|V(G_0)|$; on the other hand, it
is easy to compute that
\[
M_{k,l}=\frac{2^{k+l-1}}{k+l-1},
\]
and this matrix (up to row and column scaling, the {\it Hilbert
matrix}) has infinite rank (see e.g \cite{Cho}).

It is easy to see that in the limit, we are considering
order-preserving maps of the poset $P$ into the interval $[0,1]$;
equivalently, the limit object is the characteristic function
$W:~[0,1]^2\to[0,1]$ of the set $\{(x,y)\in[0,1]^2:~ |x-y|\ge 1/2\}$.

\section{Tools}

\subsection{Distances of functions, graphs and matrices}

For any integrable function $U:~[0,1]^2\to \R$, we define its {\it
rectangle norm} by
\begin{equation}\label{SQUAREDEF}
\|U\|_\square=\sup_{\ontop{A\subseteq [0,1]}{B\subseteq
[0,1]}}\left|\int_A\int_B U(x,y)\,dx\,dy\right|.
\end{equation}
It is easy to see that this norm could be defined by the formula
\begin{equation}\label{SQUAREALT}
\|U\|_\square= \sup_{0\le f,g \le 1} \int_0^1\int_0^1 U(x,y)f(x)g(y).
\end{equation}

The rectangle norm is related to other norms known from analysis. It
is not hard to see that
\[
\frac{1}{4} \|U\|_{\infty\to1} \le \|U\|_\square \le
\|U\|_{\infty\to1},
\]
where
\[
\|U\|_{\infty\to1} = \sup_{-1\le f,g\le 1} \int_0^1\int_0^1
U(x,y)f(x)g(y)
\]
is the $L_\infty\to L_1$ norm of the operator defined by
\[
f\mapsto \int_0^1 U(.,y)f(y)\,dy.
\]
It is also easy to see that
\[
\|U\|_\square\le \|U\|_1,
\]
where
\[
\|U\|_1=\int_0^1\int_0^1 |U(x,y)|\,dx\,dy
\]
is the $L_1$-norm of $U$ as a function.

The following lemma relates the rectangle norm and homomorphism
densities.

\begin{lemma}\label{HOMSQUARE}
Let $U,W:~[0,1]^2\to[0,1]$ be two symmetric integrable functions.
Then for every simple finite graph $F$,
\[
|t(F,U)-t(F,W)|\le |E(F)|\cdot\|U-W\|_\square.
\]
\end{lemma}

\begin{proof}
Let $V(F)=[n]$ and $E(F)=\{e_1,\dots,e_m\}$. Let $e_t=i_tj_t$. Define
$E_t=\{e_1,\dots,e_t\}$. Then
\[
t(F,U)-t(F,W)=\int_{[0,1]^n} \Bigl(\prod_{ij\in E(F)}
W(x_i,x_j)-\prod_{ij\in E(F)} U(x_i,x_j)\Bigr)\,dx\Bigr.
\]
We can write
\[
\prod_{ij\in E(F)} W(x_i,x_j)-\prod_{ij\in E(F)}
U(x_i,x_j)=\sum_{t=0}^{m-1} X_t(x_1,\dots,x_n),
\]
where
\[
X_t(x_1,\dots,x_n) = \Bigl(\prod_{ij\in E_{t-1}} W(x_i,x_j)\Bigr)
\Bigl(\prod_{ij\in E(F)\setminus E_{t}}
U(x_i,x_j)\Bigr)(W(x_{i_t},x_{j_t})- U(x_{i_t},x_{j_t})).
\]
To estimate the integral of a given term, let us integrate first the
variables $x_{i_t}$ and $x_{j_t}$; then by (\ref{SQUAREALT}),
\[
\left|\int_0^1\int_0^1 X_t(x_1,\dots,x_n)\,dx_{i_t}\,dx_{j_t} \right|
\le \|U-W\|_\square,
\]
and so
\[
|t(F,U)-t(F,W)|\le \sum_{t=0}^{m-1}\left|\int_{[0,1]^n}
X_t(x_1,\dots,x_n) \,dx \right| \le m\|U-W\|_\square
\]
as claimed.
\end{proof}

Let $G$ and $G'$ be two edge-weighted graphs on the same set $V$ of
nodes. We define their {\it  rectangular distance} by
\[
d_\square(G,G')= \frac{1}{n^2} \max_{S,T\subseteq V} \Bigl|\sum_{i\in
S,\,j\in T} (\beta_G(i,j)-\beta_{G'}(i,j))\Bigr|.
\]
Clearly, this is a metric, and the distance of any two graphs is a
real number between 0 and 1.

Finally, the {\it rectangular norm} (also called {\it cut norm} of a
matrix $A=(a_{ij})_{i,j=1}^n$ is defined by
\[
\|A\|_\square = \max_{S,T} \left| \sum_{i\in S}\sum_{j\in T}
a_{ij}\right|,
\]
where $S$ and $T$ range over all subsets of $[n]$.

These norms and distances are closely related. If $G$ and $G'$ are
two graphs on the same set of nodes, then their distance can be
expressed in terms of the associated symmetric functions $W_G$ and
$W_{G'}$, and in terms of their (weighted) adjacency matrices $A_G$
and $A_{G'}$ as
\[
d_\square(G,G')=\|W_G-W_{G'}\|_\square=\|A_G-A_{G'}\|_\square.
\]
Hence by Lemma \ref{HOMSQUARE},
\begin{equation}\label{LEFTDIST}
|t(F,G)-t(F,G')|\le |E(F)|\cdot d_\square(G,G')
\end{equation}
for any simple graph $F$.

\subsection{Szemer\'edi partitions}

A weak form of Szemer\'edi's lemma (see e.g. \cite{FK}; this weak
form is all we need) asserts that every graph $G$ can be approximated
by a weighted graph with a special structure. Let
$\PP=(V_1,\dots,V_k)$ be a partition of a finite set $V$ and let $Q$
be a symmetric $k\times k$ matrix with all entries between 0 and 1.
We define the graph $K(\PP,Q)$ as the complete graph on $V$ with a
loop at each node, in which the weight of an edge $uv$ ($u\in V_i,
v\in V_j$) is $Q_{ij}$.

\begin{lemma}[Weak form of Szemer\'edi's Lemma]\label{WEAK-SZEM}
For every $\eps>0$ there is an integer $k(\eps)>0$ such that for
every simple graph $G$ there exists a partition $\PP$ of $V(G)$ into
$k\le k(\eps)$ classes $V_1,\dots,V_k$, and a symmetric $k\times k$
matrix $Q$ with all entries between 0 and 1, such that
\[
\bigl||V_i|-|V_j|\bigr|\le 1 \qquad (1\le i,j\le k),
\]
and for every set $S\subseteq V(G)$,
\[
d_\square(G, K(\PP,Q)) \le \eps.
\]
\end{lemma}

We call the partition $\PP$ a {\it weak $\eps$-regular partition of
$G$ with density matrix $Q$.}

The best known bound $k(\eps)$ is of the order $2^{O(1/\eps^2)}$. If
the number of nodes of $G$ is less than this, then $\PP$ can be
chosen to be the partition into singletons, and $K(\PP,Q)=G$.

It is not hard to see that (at the cost of increasing $k(\eps)$) we
can impose additional conditions on the partition $\PP$. We'll need
the following condition: the partition $\PP$ refines a given
partition $\PP_0$ of $V(G)$ (the value $k(\eps)$ will also depend on
the number of classes in $\PP_0$).

It follows from the results of this paper (but it would not be hard
to prove it directly), that the above weak form of Szemer\'edi's
Lemma extends to the limit objects in the following form. A symmetric
function $U\to [0,1]^2\to[0,1]$ is called a {\it symmetric
stepfunction with $k$ steps} if there exists a partition
$[0,1]=S_1\cup\dots\cup S_k$ such that $U$ is constant on every set
$S_i\times S_j$.

\begin{corollary}\label{WEAK-SZEM-W}
For every $\eps>0$ there is an integer $k(\eps)>0$ such that for
every symmetric function $W:~[0,1]^2\to[0,1]$ there exists a
symmetric stepfunction $U:~[0,1]^2\to[0,1]$ with $k$ steps such
\[
\|W-U\|_\square\le\eps.
\]
\end{corollary}

(It would be interesting to find a similar form of the full-strength
Szemer\'edi Lemma.)

\section{Proofs}\label{PROOFS}

\subsection{Proof of Lemma \ref{HOM-INJ}}
We have trivially
\[
\hom(F,G)\ge\inj(F,G),
\]
and so
\[
t(F,G)=\frac{\hom(F,G)}{n^k} \ge
\frac{\inj(F,G)}{n^k}=t_0(F,G)\frac{(n)_k}{n^k}.
\]
Here
\[
\frac{(n)_k}{n^k} = \prod_{i=1}^{k-1} \left(1-\frac{i}{n}\right)\ge
1-\binom{k}{2}\frac{1}{n},
\]
and so
\[
t(F,G) \ge t_0(F,G)\left(1-\binom{k}{2}\frac{1}{n}\right) \ge
t_0(F,G) -\binom{k}{2}\frac{1}{n}.
\]
On the other hand, we have by the beginning of inclusion-exclusion,
\[
\inj(F,G)\ge \hom(F,G)-\sum_{F'} \hom(F',G),
\]
where the summation ranges over all graphs $F'$ arising from $F$ by
identifying two of the nodes. The number of such graphs is
$\binom{k}{2}$. Hence
\begin{align*}
t_0(F,G)&=\frac{\inj(F,G)}{(n)_k} \ge \frac{\inj(F,G)}{n^k} \ge
\frac{\hom(F,G)}{n^k} - \sum_{F'} \frac{\hom(F',G)}{n^k}\\ &=
t(F,G)-\frac{1}{n}\sum_{F'} t(F',G)\ge
t(F,G)-\binom{k}{2}\frac{1}{n}.
\end{align*}
This completes the proof.

\subsection{Proof of Lemma \ref{EXPINJ}}

Consider any injective map $\phi:~V(F)\to V(\Ge)$. For a fixed choice
of $X_1,\dots,X_n$, the events $\phi(i)\phi(j)\in E(\Ge)$ are
independent for different edges $ij$ of $F$, and so the probability
that $\phi$ is a homomorphism is
\[
\prod_{ij\in E(F)} W(X_{\phi(i)},X_{\phi(j)})
\]
Now choosing $X_1,\dots,X_n$ at random, we get that the probability
that $\phi$ is a homomorphism is
\[
\E\Bigl(\prod_{ij\in E(F)} W(X_{\phi(i)},X_{\phi(j)})\Bigr)=t(F,W).
\]
Summing over all injective maps $\phi$, we get (a). By
(\ref{HOM-INJ}), we get (b).

Finally, we estimate the variance of $t(F,\Ge)$. Let $F_2$ denote the
disjoint union of $2$ copies of $F$. Then
\[
t(F_2,\Ge)=t(F,\Ge)^2, \qquad\text{and}\qquad t(F_2,W)=t(F,W)^2.
\]
Let $R=|V(F)|^2/n$. By Lemma \ref{HOM-INJ}
\[
\E(t(F,G)^2)=\E(t(F_2,G))\le \E(t_0(F_2,\Ge)+
2R)=t(F_2,W)+2R=t(F,W)^2+2R,
\]
and
\[
\E(t(F,G))^2\ge \E(t_0(F,G)-R/2)^2 \ge
\E(t_0(F,\Ge))^2-R=t_0(F,W)^2-R.
\]
Hence
\[
\Var(t(F,\Ge))= \E(t(F,G)^2)-\E(t(F,\Ge))^2\le
3R=\frac{3}{n}|V(F)|^2.
\]

\subsection{Proof of Theorem \ref{CONC}}

The idea of the proof is to form a martingale as follows. In the
$m$-th step ($m=1,\dots,n$), we generate $X_m\in[0,1]$, and the edges
of $\Ge$ connecting the new node to previously generated nodes. The
probability that a random injection of $V(F)$ into $V(\Ge)$ is a
homomorphism (conditioning on the part of $\Ge$ we already generated)
is a martingale. We are going to apply Azuma's inequality to this
martingale.

To be precise: For every injective map $\phi:~[k]\to[m]$, let
$A_\phi$ denote the event that $\phi$ is a homomorphism from $F$ to
the random graph $\Ge$. Let $\Ge_m$ denote the subgraph of $\Ge$
induced by nodes $1,\dots,m$. Define
\[
B_m=\frac{1}{(n)_k}\sum_{\phi} \Pr(A_\phi~|~\Ge_m).
\]
Clearly the sequence $(B_0,B_1,\dots)$ is a martingale. Furthermore,
\[
\Pr(A_\phi)=t(F,W),
\]
and
\[
\Pr(A_\phi~|~\Ge_n)=
  \begin{cases}
    1 & \text{if $\phi$ is a homomorphism from $F$ to $\Ge$}, \\
    0 & \text{otherwise}.
  \end{cases}
\]
Thus
\[
B_0=\sum_\phi \Pr(A_\phi)=t(F,W),
\]
and
\[
B_n=\frac{1}{(n)_k}\inj(F,\Ge)=t_0(F,\Ge).
\]

Next we estimate $|B_m-B_{m-1}|$:
\begin{align*}
|B_m-B_{m-1}| &=\frac{1}{(n)_k} \left|\sum_{\phi}
(\Pr(A_\phi~|~\Ge_m)-\Pr(A_\phi~|~\Ge_{m-1}))\right|\\
&\le \frac{1}{(n)_k}\sum_{\phi}
\Bigl|\Pr(A_\phi~|~\Ge_m)-\Pr(A_\phi~|~\Ge_{m-1})\Bigr|.
\end{align*}
In this sum, every term for which $m$ is not in the range of $\phi$
is 0, and the other terms are at most 1. The number of terms of the
latter kind is $k(n-1)_{k-1}$, and so
\[
|B_m-B_{m-1}| \le \frac{k(n-1)_{k-1}}{(n)_k}=\frac{k}{n}.
\]

Thus we can invoke Azuma's Inequality:
\[
\Pr\Bigl(B_n-B_0>\eps\Bigr) \le
\exp\left(-\frac{\eps^2}{2n(k/n)^2}\right)=
\exp\left(-\frac{\eps^2}{2k^2}n\right),
\]
and similarly
\[
\Pr\Bigl(B_n-B_0>\eps(n)_k\Bigr) \le
\exp\left(-\frac{\eps^2}{2k^2}n\right).
\]
Hence
\[
\Pr\Bigl(|B_n-(n)_kt(F,W)|>\eps(n)_k\Bigr) \le
2\exp\left(-\frac{\eps^2}{2k^2}n\right).
\]
This proves (\ref{INJFH}).

To get (\ref{TFH}), we use Lemma \ref{HOM-INJ}. We may assume that
$n>k^2/\eps$ (else the inequality is trivial). Then
\[
|t_0(F,\Ge)-t(F,\Ge)|\le \frac{\eps}{3},
\]
and similarly
\[
|t_0(F,\Ge)-t(F,H)|\le \frac{\eps}{3},
\]
so (\ref{TFH}) follows by applying (\ref{INJFH}) with $\eps/3$ in
place of $\eps$.

\subsection{Proof of Theorem \ref{MAIN}: (a)$\Rightarrow$(b)}

Let $(G_n)$ be a convergent graph sequence and
\[
f(F)=\lim_{n\to\infty} t(F,G_n)
\]
for every $n$. We want to construct a function $W:~[0,1]^2\to[0,1]$
such that $f=t(.,W)$.

We start with constructing a subsequence of $(G_n)$ whose members
have well-behaved Szemer\'edi partitions.

\begin{lemma}\label{SUBSEQ}
Every graph sequence $(G_n:~n=1,2,\dots)$ has a subsequence
$(G'_m:~m=1,2,\dots)$ for which there exists a sequence
$(k_m:~m=1,2,\dots)$ of integers and a sequence $(Q_m:~m=1,2,\dots)$
of matrices with the following properties.

\smallskip

{\rm (i)} $Q_m$ is a $k_m\times k_m$ symmetric matrix, all entries of
which are between $0$ and $1$.

\smallskip

{\rm (ii)} If $i<j$, then $k_i\mid k_j$, and the matrix $Q_i$ is
obtained from the matrix $Q_j$ by partitioning its rows and columns
into $k_i$ consecutive blocks of size $k_j/k_i$, and replacing each
block by a single entry, which is the average of the entries in the
block.

\smallskip

{\rm (iii)} For all $j<m$, $G'_m$ has a weakly $(1/m)$-regular
partition $\PP_{m,j}$ with density matrix $Q_{m,j}$ such that
\begin{equation}\label{QDIST}
\|Q_{m,j}-Q_j\|_\square<1/j,
\end{equation}
and for $1\le i<j\le m$, $\PP_{m,j}$ is a refinement of $\PP_{m,i}$.
\end{lemma}

\begin{proof}
For every integer $m\ge 1$, we construct a subsequence
$(G_n^m)$ so that all graphs in the subsequence have a weakly
$(1/m)$-regular partition $\PP_{n,m}$ into the same number $k_m$ of
classes and with almost the same density matrix.

The first sequence $(G_n^1)$ is selected from $(G_n)$ so that the
edge density of $G_n^1$ converges to a fix constant $c$ between $0$
and $1$ if $n$ tends to infinity. Furthermore, for every graph
$G_n^1$ let $\PP_{n,1}=\{V(G_n^1)\}$ be the $1$-block partition and
let $Q_{n,1}$ be the $1$ by $1$ matrix containing the edge density of
$G_n^1$. We set $k_1=1$.

Suppose that for some integer $m>0$, we have constructed the sequence
$(G_n^m)$. For every graph $G_n^m$, consider a weakly
$1/(m+1)$-regular partition $\PP_{n,m+1}=\{V_1,\dots, V_{K_n}\}$ of
$G^m_n$ with density matrix $Q_{n,{m+1}}$. We may choose this
partition so that it refines the previous partition $\PP_{n,m}$, and
each class of $\PP_{n,m}$ is split into the same number of classes
$r_{n,m}$; the number $K_n$ of classes remains bounded, and so the
numbers $r_{n,m}$ also remain bounded ($m$ is fixed, $n\to\infty$).
So we can thin the sequence so that all remaining graphs have the
same $r_{m}=r_{n,m}$. We set $k_{m+1}=k_mr_m$. Furthermore, we can
select a subsequence so that the density matrices $Q_{n,m+1}$
converge to a fixed matrix $Q_{m+1}$ if $n$ tends to infinity.
Finally we drop all the elements $G_i^{m+1}$ from the remaining
sequence for which $\|Q_{i,m+1}-Q_{m+1}\|_\square >1/(2m+2)$. By
renumbering the indices, we obtain the subsequence $G_n^{m+1}$.

Let $G'_m=G_1^m$. For every $1\le j<m$, the graph $G'_m$ has a tower
of partitions $\PP_{m,j}$ ($j=1,\dots,m$) so that $\PP_{m,j}$ has
$k_j$ almost equal classes, and has a density matrix $Q_{m,j}$ such
that $\lim_{m\to\infty}Q_{m,j}=Q_j$ and $\|Q_{m,j}-Q_{m',j}\|_\square
<1/j$ for all $m,m'>j$.

For a fixed graph $G'_m$, the partitions
$\PP_{m,1},\PP_{m,2},\dots,\PP_{m,m}$ are successively refinements
of each other. We may assume that the classes are labeled so that
the $i$-th class of $\PP_{m,j+1}$ is the union of consecutive
classes $(i-1)r_m+1,\dots,ir_m$. Let $\widehat{Q}_{m,j}$ ($1\le
j\le m$) be the $k_j\times k_j$ matrix obtained from $Q_m$ by
partitioning its rows and columns into $k_j$ consecutive blocks of
size $k_m/k_j$, and replacing each block by a single entry, which
is the average of the corresponding entries of $Q_m$. Using that
$\lim_{n\to\infty} |V(G_n)|=\infty$ and that all the sets in
$\PP_{m,j}$ have almost the same size one gets that
$\widehat{Q}_{m,j}=Q_j$.

Thus we constructed a sequence $(G'_1,G'_2,\dots)$ of graphs, an
increasing sequence $(k_1,k_2,\dots)$ of positive integers, and a
sequence $(Q_1,Q_2\dots)$ of matrices. We claim that these sequences
satisfy the properties required in the Lemma. (i), (ii) and the
second assertion of (iii) are trivial by construction. The first
assertion in (iii) follows on noticing that $Q_j$ is the limit of
matrices $Q_{m,j}$, and $\|Q_{m,j},Q_{m',j}\|_\square<1/j$ for all
$j\leq m \leq m'$.
\end{proof}

\begin{lemma}\label{W-EXISTS}
Let $(k_m)$ be a sequence of positive integers and $(Q_m)$, a
sequence of matrices satisfying {\rm (i)} and {\rm (ii)} in Lemma
\ref{SUBSEQ}. Then there exists a symmetric measurable function
$W:~[0,1]^2\to [0,1]$ such that

{\rm (a)} $W_{Q_m}\to W~ (m\to\infty)$ almost everywhere;

\smallskip

{\rm (b)} for all $m$ and $1\le i,j\le k_m$
\[
(Q_M)_{ij}=k_m^2\int_{(i-1)/k_m}^{i/k_m}\int_{(j-1)/k_m}^{j/k_m}
W(x,y)\,dx\,dy.
\]
\end{lemma}

\begin{proof}
Define a map $\phi_m:~[0,1)\to[k_m]$ by mapping the interval
$[(i-1)/k_m,i/k_m)$ to $i$.

\begin{claim}\label{MARTINGAL}
Let $X$ and $Y$ be two uniformly distributed random elements of
$[0,1]$. Then the sequence $Z_m=(Q_m)_{\phi_m(X),\phi_m(Y)}$,
$i=1,2,\dots$ is a martingale.
\end{claim}

We want to show that
\begin{equation}\label{MARTIN-1}
\E(Z_{m+1}\mid Z_1,\dots,Z_m)=Z_m.
\end{equation}
In fact, we show that
\begin{equation}\label{MARTIN-2}
\E(Z_{m+1}\mid \phi_1(X),\phi_1(Y),\dots,\phi_m(X),\phi_m(Y))=Z_m.
\end{equation}

Since $\phi_m(X)$ and $\phi_m(Y)$ determine $\phi_i(X)$ and
$\phi_i(Y)$ for $i<m$, it suffices to show that
\begin{equation}\label{MARTIN-3}
\E(Z_{m+1}\mid \phi_m(X)=a,\phi_m(Y)=b)=(Q_m)_{a,b}.
\end{equation}
The condition $\phi_m(X)=a$ and $\phi_m(Y)=b$ force $X$ to be uniform
in the interval $[(a-1)/k_m,a_{k_m}]$, and so $\phi_{m+1}(X)$ is a
uniform integer in the interval
$[(k_{m+1}/k_m)(a-1),(k_{m+1}/k_m)a]$. Similarly, $\phi_{m+1}(Y)$ is
a uniform integer in the interval
$[(k_{m+1}/k_m)(b-1),(k_{m+1}/k_m)b]$. So
$(Q_{m+1})_{\phi_{m+1}(X),\phi_{m+1}(Y)}$ is a uniformly distributed
entry of the submatrix formed by these rows and columns. By condition
(ii), the average of these matrix entries is exactly $(Q_m)_{a,b}$.
This proves the claim.

Since $Z_m$ is also bounded, we can invoke the Martingale Convergence
Theorem, and conclude that $\lim_{m\to\infty} Z_m$ exists with
probability 1. This means that
\[
W(x,y)=\lim_{m\to\infty} (Q_m)_{\phi_m(x),\phi_m(y)}
\]
exists for almost all pairs $(x,y)$, $0\le x,y\le 1$. Let us define
$W(x,y)=0$ whenever the limit does not exist.

It is trivial that $W$ is symmetric, $0\le W\le 1$, and $W$ satisfies
condition (a) in the Lemma. Furthermore,
\begin{align*}
\int_{(i-1)/k_m}^{i/k_m}&\int_{(j-1)/k_m}^{j/k_m} W(x,y)\,dx\,dy =
\int_{(i-1)/k_m}^{i/k_m}\int_{(j-1)/k_m}^{j/k_m} \lim_{n\to\infty}
(Q_n)_{\phi_n(x),\phi_n(y)}\,dx\,dy\\ &=\lim_{n\to\infty}
\int_{(i-1)/k_m}^{i/k_m}\int_{(j-1)/k_m}^{j/k_m}
(Q_n)_{\phi_n(x),\phi_n(y)}\,dx\,dy \\
&=\frac{1}{k_m^2}\Bigl(\lim_{n\to\infty}\frac{k_m^2}{k_n^2}
\sum_{a=(i-1)(k_n/k_m)+1}^{i(k_n/k_m)}
\sum_{b=(j-1)(k_n/k_m)+1}^{j(k_n/k_m)}(Q_n)_{a,b}\Bigr)\\
&=\frac{1}{k_m^2}(Q_m)_{i,j}
\end{align*}
(the last equation follows from assumption (ii)). This proves (b).
\end{proof}

Now it easy to conclude the proof of the necessity of the condition
in Theorem \ref{MAIN}. Let us apply Lemma \ref{SUBSEQ} to the given
convergent graph sequence. The sequence $(G'_1,G'_2,\dots)$ of graphs
it gives is a subsequence of the original sequence, so it is
convergent and defines the same limit parameter $f$. The lemma also
gives a sequence of integers and a sequence of matrices satisfying
(i) and (ii). We can use Lemma \ref{MARTINGAL} to construct an
integrable function $W:~[0,1]^2\to[0,1]$ with properties (a), (b) and
(c).

It remains to show that $f=t(.,W)$. For $1\le j\le m$, let
$G^*_{m,j}=G(\PP_{m,j},Q_{m,j})$ and $G_{m,j}^{**}=G(\PP_{m,j},Q_j)$.
Then
\begin{equation}\label{GM-GMJ}
d(G'_m,G^*_{m,j})\le \frac{1}{j}
\end{equation}
(since $\PP_{m,j}$ is a weakly $(1/j)$-regular partition of $G'_m$),
and
\begin{equation}\label{GMJ-GMJ}
d(G^*_{m,j},G^{**}_{m,j})\le \frac{1}{j}
\end{equation}
by (\ref{QDIST}).

Let $W_{m,j}=W_{G_{m,j}^{**}}$ and $W_j=W_{Q_j}$. Clearly
\[
t(F,W_{mj})=t(F,G_{m,j}^{**}).
\]
Furthermore,
\begin{equation}\label{WMJ-WJ}
W_{m,j}\to W_j\qquad (m\to\infty)
\end{equation}
almost everywhere. Indeed, the functions $W_{m,j}$ and $W_j$ differ
only if the classes in $\PP_{m,j}$ are not all equal; but even in
this case, if $W_{m,j}(x,y)\not=W_j(x,y)$ then either $x$ or $y$ must
be closer to one of the numbers $a/k_j$ than $1/|V(G_m')|$. Finally,
we have
\begin{equation}\label{WJ-W}
W_j\to W\qquad (j\to\infty)
\end{equation}
almost everywhere.

Now let $\eps>0$, and choose a positive integer $m_0$ so that for
$m>m_0$, we have
\[
|t(F,G'_m)-f(F)|<\frac{\eps}{4}.
\]
By (\ref{WJ-W}), we can choose a positive integer $j$ so that
\[
|t(F,W)-t(F,W_j)|<\frac{\eps}{4}.
\]
We may also assume that $j>8/\eps$ and $j>m_0$. By (\ref{WMJ-WJ}), we
can choose an $m>j$ so that
\[
|t(F,W_{m,j})-t(F,W_j)|\le\frac{\eps}{4}.
\]
By Lemma \ref{LEFTDIST} and inequalities (\ref{GM-GMJ}) and
(\ref{GMJ-GMJ}), we have
\[
|t(F,G'_m)-t(F,G^{**}_{m,j})|\le\frac{\eps}{4}.
\]
Combining these inequalities, we get that
\[
|f(F)-t(F,W)|\le\eps,
\]
which completes the proof of (a)$\Rightarrow(b)$.

\subsection{Proof of Theorem \ref{MAIN}: (b)$\Rightarrow$(c)}

Let $f=t(.,W)$. It is obvious that $f$ is normalized and
multiplicative.

To prove that $f$ is reflection positive, consider any finite set
$F_1,\dots,F_m$ of $k$-labeled graphs, and real numbers
$y_1,\dots,y_m$. We want to prove that
\[
\sum_{p,q=1}^m f(F_pF_q) y_py_q \ge 0.
\]
For every $k$-labeled graph $F$ with node set $[n]$, let $F'$ denote
the subgraph of $F$ induced by the labeled nodes, and $F''$ denote
the graph obtained from $F$ by deleting the edges spanned by the
labeled nodes. Define
\[
\tau(F,x_1,\dots,x_k)=\int_{[0,1]^{n-k}} \prod_{ij\in E(F'')}
W(x_i,x_j)\,dx_{k+1}\dots dx_n,
\]
and for every graph $F$ with $V(F)=[k]$,
\[
W(F,x_1,\dots,x_k)=\prod_{ij\in E(F)} W(x_i,x_j).
\]
Then
\begin{align*}
\sum_{p,q=1}^m y_py_q f(F_pF_q) = \int_{[0,1]^k}& \sum_{p,q=1}^m
y_p y_q\tau(F_p,x_1,\dots,x_k)\tau(F_q,x_1,\dots,x_k)\\
&W(F'_p\cup F'_q,x_1,\dots,x_k)\,dx_1\dots dx_k.
\end{align*}
We prove that the integrand is nonnegative for every $x_1,\dots,x_k$:
\[
\sum_{p,q=1}^m \overline{y}_p\overline{y}_qW(F'_p\cup F'_q)\ge 0,
\]
where $\overline{y}_p=y_p\tau(F_p,x_1,\dots,x_k)$, and the $x_i$ are
suppressed for clarity). Let
\[
\hat W(F)= \prod_{ij\in E(F)} W(x_i,x_j)\prod_{ij\in E(\overline{F})}
(1-W(x_i,x_j)),
\]
then clearly $W(\overline{F})\ge 0$, and for every $F\in\FF_k$,
\[
W(F)=\sum_{H\supset F} \hat W(H)
\]
(where the summation extends over all $H\in\FF_k$ containing $F$ as a
subgraph). Thus
\begin{align*}
\sum_{p,q=1}^m &\overline{y}_p\overline{y}_qW(F'_p\cup F'_q) =
\sum_{p,q=1}^m \overline{y}_p\overline{y}_q \sum_{H\supset F'_p\cup
F'_q} \hat W(H)\\
&= \sum_{H\in\FF_k}\hat W(H) \sum_{p,q: F_p,F_q\subseteq H}
\overline{y}_p\overline{y}_q =\sum_{H\in\FF_k}\hat W(H)
\left(\sum_{p: F_p\subseteq H} \overline{y}_p\right)^2 \ge 0.
\end{align*}

\subsection{Proof of Theorem \ref{MAIN}: (c)$\Rightarrow$(d)}

This is trivial, since $M_0(k,f)$ is a symmetric submatrix of
$M(k,f)$.

\subsection{Proof of Theorem \ref{MAIN}: (d)$\Rightarrow$(e)}

The proof of this implication uses the fact that the entry in the
$(F_1,F_2)$-position in $M_0(k,f)$ depends on the union of $F_1$ and
$F_2$ only. The Lindstr\"om--Wilf Formula gives a nice
diagonalization of such matrices  as follows. Let $\FF_k$ denote the
set of all graphs with nodes $[k]$. Let $Z$ denote the
$\FF_k\times\FF_k$ matrix defined by
\begin{equation}\label{LWZ}
Z_{F_1,F_2}=
  \begin{cases}
    1 & \text{if $F_1\subseteq F_2$}, \\
    0 & \text{otherwise}.
  \end{cases}
\end{equation}
Let $D$ be the diagonal matrix
\begin{equation}\label{LWD}
D_{F_1,F_2}=
  \begin{cases}
    f^\dagger(F_1) & \text{if $F_1=F_2$}, \\
    0 & \text{otherwise}.
  \end{cases}
\end{equation}
Then
\begin{equation}\label{LWMAIN}
B=Z\T D Z.
\end{equation}
This implies that $M_0(k,f)$ is positive semidefinite if and only
if $f^{\dagger}\ge 0$ for all graphs with $k$ nodes.

\subsection{Proof of Theorem \ref{MAIN}: (e)$\Rightarrow$(a)}

Let $f$ be a normalized and multiplicative graph parameter such that
$f^{\dagger}\ge 0$. Fix any $k\ge 1$. As a first step we construct a
random variable $\Ge_k$, whose values are graphs with $k$ labeled
nodes: Let $\Ge_k=F$ with probability $f^{\dagger}(F)$. Since
$f^\dagger\ge 0$ by hypothesis and
\[
\sum_F f^{\dagger}(F) = f(\overline{K_k})=1
\]
(where the summation extends over all graphs $F$ with $V(F)=[k]$),
this is well defined. It is also clear that this distribution does
not depend on the labeling of the nodes.

Next we show that for every graph $F$ with $k$ nodes,
\begin{equation}\label{FIXN}
f(F)=\E(t_0(F,\Ge_k)).
\end{equation}
Indeed, we have
\begin{align*}
\E(t_0(F,\Ge_k)) &= \sum_{F'\supseteq F} \E(t_1(F',\Ge_k)) \\
&=\sum_{F'\supseteq F} \Pr(F'=\Ge_k)= \sum_{F'\supseteq
F}f^{\dagger}(F') =f(F).
\end{align*}

We claim that for every graph with $k\le n$ nodes,
\begin{equation}\label{EXP1}
f(F)=\E(t_0(F,\Ge_n)).
\end{equation}
Indeed, add $n-k$ isolated nodes to $F$ to get a graph $F'$ with $n$
nodes. Then $f(F')=f(F)$ by multiplicativity and $f(K_1)=1$, while
$t_0(F',G)=t_0(F,G)$ for every graph $G$. Thus
\[
f(F)=f(F')=\E(t_0(F',\Ge_n))=\E(t_0(F,\Ge_n)).
\]

We need a bound on the variance of $t_0(F',\Ge_n)$: By (\ref{EXP1}),
\[
\Var(t_0(F,\Ge_n))=\E(t_0(F,\Ge_n)^2)-(\E(t_0(F,\Ge_n))^2.
\]
Here
\[
(\E(t_0(F,\Ge_n))^2= f(F)^2=f(FF)=\E(t_0(FF,\Ge_n)
\]
(by multiplicativity), so
\[
\Var(t_0(F,\Ge_n))=\E\Bigl(t_0(F,\Ge_n)^2-t_0(FF,\Ge_n)\Bigr).
\]
Now for any graph $G$,
\[
t(F,G)^2=t(FF,G),
\]
and so
\begin{align*}
|t_0(F,\Ge_n)^2-t_0(FF,\Ge_n)| &\le
|t(F,\Ge_n)^2-t_0(F,\Ge_n)^2|+|t(FF,\Ge_n)-t_0(FF,\Ge_n)| \\
&\le 2|t(F,\Ge_n)-t_0(F,\Ge_n)|+|t(FF,\Ge_n)-t_0(FF,\Ge_n)|\\
& \le2\binom{k}{2}\frac{1}{n}+\binom{2k}{2}\frac{1}{n}
<\frac{3k^2}{n}.
\end{align*}
Thus for every graph $F$ with $k\le n$ nodes,
\begin{equation}\label{VAR1}
\Var(t_0(F,\Ge_n))\le \frac{3k^2}{n}.
\end{equation}

By Chebyshev's Inequality,
\[
\Pr(|t_0(F,\Ge_n)-f(F)|>\eps) <\frac{3k^2}{\eps^2 n}.
\]
It follows by the Borell-Cantelli Lemma that if we take (say) the
graph sequence $(\Ge_{n^2}:~n=1,2,\dots)$, then with probability 1,
\[
t_0(F,\Ge_{n^2})\to f(F)\quad (n\to\infty).
\]
Since there are only a countable number of graphs $F$, this
convergence holds with probability 1 for every $F$. So we see that
\[
f(.)=\lim_{n\to\infty} t(.,\Ge_{n^2})
\]
for almost all choices of the sequence $(\Ge_{n^2})$. This completes
the proof of Theorem \ref{MAIN}.

\medskip

\noindent{\bf Remarks. 1.} There are alternatives for certain parts
of the proof. Instead of verifying (b)$\Rightarrow$(c) directly, we
could argue that (b)$\Rightarrow$(a) (which follows e.g. from
Corollary \ref{W-CONVERGE}), and then that (a)$\Rightarrow$(c) (which
follows from the characterization of homomorphism functions in
\cite{FLS}).

\medskip

{\bf 2.} Equation (\ref{EXP1}), satisfied by the random graph
$\Ge_n$, is the same as equation (a) in Lemma \ref{EXPINJ}, satisfied
by the random graph $\Ge(n,W)$. It is not hard to see that this
equation uniquely determines the distribution on $n$-node graphs, and
hence $\Ge_n$ and $\Ge(n,W)$ have the same distribution.

\medskip

{\bf 3.} The construction of the limit object $W$ in the proof shows
that every convergent graph sequence $(G_n)$ has a subsequence
$(G_n')$ such that (with an appropriate labeling of the nodes)
\[
\|W_{G_n'}-W\|_\square \to 0.
\]

\subsection{Proof of Theorem \ref{GEN-RAND}}

It is trivial that $\Ge(n,W)$ satisfies (a), (b) and (c). Conversely,
suppose that $\Ge_n$ is a graph model with properties (i), (ii) and
(iii). Define a graph parameter $f$ by
\[
f(F)=\Pr(F\subseteq \Ge_k),
\]
where $V(F)=[k]$. By condition (i), $f$ is independent of the
labeling of the nodes, so it is indeed a graph parameter.

We claim that $f\in\TT$, by verifying (e) in Theorem \ref{MAIN}. It
is trivial that $f$ is normalized. Multiplicativity is an easy
consequence of (iii) and (ii): Let $V(F_1)=[k]$ and
$V(F_2)=\{k+1,\dots,k+l\}$, then
\begin{align*}
f(F_1\cup F_2)&=\Pr(F_1\cup F_2\subseteq
\Ge_{k+l})=\Pr(F_1\subseteq \Ge_{k+l})\Pr(F_2\subseteq
\Ge_{k+l})\\&=\Pr(F_1\subseteq \Ge_{k})\Pr(F_2\subseteq
\Ge_{l})=f(F_1)f(F_2).
\end{align*}

Next we show that
\begin{equation}\label{FETG}
f(F)=\E(t_0(F,\Ge_n))
\end{equation}
for every $n\ge k$. For $n=k$, this is obvious. Let $n>k$, and let us
add $n-k$ isolated nodes to $F$ to get graph $F'$. Then, using
condition (ii),
\[
f(F)=\Pr(F\subseteq \Ge_k)=\Pr(F'\subseteq
\Ge_n)=\E(t_0(F',\Ge_n))=\E(t_0(F,\Ge_n))
\]
as claimed.

Applying (\ref{FETG}) we obtain
\[
f^\dagger(F) = \E(t_0^\dagger(F,\Ge_n)) = \E(t_1(F,\Ge)) \ge 0.
\]
This proves that $f\in \TT$.

By Theorem \ref{MAIN}, there exists a symmetric measurable function
$W:~[0,1]^2\to[0,1]$ such that
\[
f(F)=t(F,W)
\]
for every graph $F$. By Lemma \ref{EXPINJ}(a) and equation
(\ref{FETG}), we have
\[
\E(t_0(F,\Ge_k))=\E(t_0(F,\Ge(k,W)).
\]
Applying the ``$\dagger$'' operator, this implies that
\[
\E(t_1(F,\Ge_k))=\E(t_1(F,\Ge(k,W)),
\]
and hence
\[
\Pr(F=\Ge_k)=\Pr(F=\Ge(k,W)),
\]
which proves that $\Ge_k$ and $\Ge(k,W)$ have the same distribution.

\section{Concluding remarks}

We mention some results and problems related to our work. Details
(exact formulations, results and conjectures) will be discussed
elsewhere.

\subsection{Uniqueness}

The limit function of a graph sequence is ``essentially unique''. In
other words, if two functions $U,W:~[0,1]^2\to[0,1]$ are such that
\[
t(F,U)=t(F,W)
\]
for every simple graph $F$, then $U$ and $W$ are ``essentially the
same''. Unfortunately, it is nontrivial to characterize what
``essentially the same'' means; for example, $U$ could be obtained
from $W$ by applying the same measure-reserving permutation in both
coordinates.

\subsection{Weighted graphs and multiple edges}\label{EXTENSIONS}

It seems to be quite straightforward to extend our results to the
case when the graphs $G_n$ can have multiple edges, or more
generally, edge-weights (not restricted to $[0,1]$): We simply have
to drop the bounds on the limit function $W$. However, several
technical issues arise concerning integrability conditions and the
applicability of the Martingale Theorem.

Allowing multiple edges in the ``sample graphs'' $F$ leads to a more
complicated question. Assume that we consider a sequence of simple
graphs $G_n$. If we define $F'$ as the underlying simple graph of a
multigraph $F$, then
\[
\hom(F',G_n)=\hom(F,G_n)\qquad \text{and}\qquad t(F',G_n)=t(F,G_n),
\]
but
\[
t(F',W)<t(F,W)
\]
if $W$ is a function that is strictly between 0 and 1. Since, as we
have remarked, the limit function $W$ is essentially unique, the
formula for $t(F,W)$ does not define the limit of $t(F,G_n)$
correctly.

\subsection{Extremal graph theory}

There are many results in graph theory, especially in extremal graph
theory, that can be formulated as inequalities between the numbers
$t(F,G)$ for a fixed $G$ and various graphs $F$. For example,
Goodman's theorem relating the number of edges to the number of
triangles can be stated as
\[
t(K_3,G)\ge t(K_2,G)(2t(K_2,G)-1).
\]
This inequality is equivalent to saying that for every graph
parameter $t\in\TT$,
\begin{equation}\label{TURAN}
t(K_3)\ge t(K_2)(2t(K_2)-1).
\end{equation}
By Theorem \ref{MAIN}, such an inequality must be a consequence of
reflection positivity, multiplicativity, and the trivial condition
that $t$ is normalized. In fact, (\ref{TURAN}) can be easily derived
from these conditions (this is left to the reader as an exercise).
Many other results in extremal graph theory follow in a similar way.

\section*{Acknowledgement}

We are grateful to Jeong-Han Kim for his kind advice on Azuma's
inequality.


\begin{thebibliography}{99}

\bibitem{BS}
I.~Benjamini and O.~Schramm: Recurrence of Distributional Limits of
Finite Planar Graphs, {\it Electronic J. Probab.} {\bf 6} (2001),
paper no. 23, 1--13.

\bibitem{Cho}
M.-D.~Choi: Tricks or Treats with the Hilbert Matrix, {\it Amer.
Math. Monthly} {\bf 90} (1983), 301--312.

\bibitem{CGW}
F.~Chung, R.L.~Graham and R.M.~Wilson: Quasi-random graphs, {\it
Combinatorica} {\bf 9} (1989), 345--362.

\bibitem{ELS}
P.~Erd\H{o}s, L.~Lov\'asz, J.~Spencer: Strong independence of
graphcopy functions, in: {\it Graph Theory and Related Topics},
Academic Press, 165-172.


\bibitem{FLS}
M.~Freedman, L.~Lov\'asz, A.~Schrijver: Reflection positivity, rank
connectivity, and homomorphism of graphs (MSR Tech Report \#
MSR-TR-2004-41)

\url{ftp://ftp.research.microsoft.com/pub/tr/TR-2004-41.pdf}

\bibitem{FK}
A.~Frieze and R.~Kannan: Quick approximation to matrices and
applications, {\it Combinatorica} {\bf 19}, 175--220.

\bibitem{Gow}
W.T.~Gowers: Lower bounds of tower type for Szemer\'edi's Uniformity
Lemma, {\it Geom. Func. Anal.} {\bf 7} (1997), 322--337.

\bibitem{LL}
L.~Lov\'asz: Operations with structures, {\it Acta Math.\ Hung.} {\bf
18}, 321-328.

\bibitem{LS}
L.~Lov\'asz and V.T.~S\'os: Generalized quasirandom graphs,
manuscript.

\bibitem{Ly}
R.~Lyons: Asymptotic enumeration of spanning trees {\it Combin.
Probab. Comput.}, to appear.

\bibitem{Tho}
A.~Thomason: Pseudorandom graphs, in: {\it Random graphs '85}
North-Holland Math. Stud. {\bf 144}, North-Holland, Amsterdam, 1987,
307--331.

\end{thebibliography}
\end{document}